\documentclass{conm-p-l}
\usepackage{amssymb,latexsym,amsmath,amscd,yhmath,xypic}
\newtheorem{thm}{Theorem}[section]
\newtheorem{prop}[thm]{Proposition}

\newtheorem{lemma}[thm]{Lemma}

\theoremstyle{definition}
\newtheorem{defn}[thm]{Definition}

\newtheorem{quest}[thm]{Question}
\newtheorem{ex}[thm]{Example}

\theoremstyle{remark}
\newtheorem{rmk}[thm]{Remark}

\numberwithin{equation}{section}

\newcommand{\lra}{\longrightarrow}
\newcommand{\PP}{\mathbb{P}}

\newcommand{\ZZ}{\mathbb{Z}}

\newcommand{\bin}[2]{ {{#1} \choose {#2}} }

\newcommand{\HH}{\mathbb{H}}
\newcommand{\I}{\mathcal{I}}

\newcommand{\MM}{\mathfrak{m}}

\newcommand{\FF}{\mathbb{F}}
\newcommand{\GG}{\mathbb{G}}

\newcommand{\F}{\mathcal{F}}

\newcommand{\mmin}{{\rm min}}
\newcommand{\hgt}{{\rm ht}}

\title{Lifting the determinantal property}

\author{Elisa Gorla}

\address{Institut f\"ur Mathematik
\\ Universit\"at Z\"urich, \hfil\break\indent Winterthurerstrasse
190, CH-8057 Z\"urich, Switzerland}

\email{elisa.gorla@math.unizh.ch}

\thanks{I am grateful to J.~Migliore and A.~Conca for useful
discussions. Part of the research in this paper was carried out
while the author was a guest at the Max Planck Institut f\"ur 
Mathematik in Bonn. The computer algebra system CoCoA \cite{cocoa} was
used for some of the computations. The author was partially supported
by the Swiss National Science Foundation under grant no. 107887.}

\begin{document}

\maketitle

\begin{abstract}
In this note we study standard and in particular good determinantal
schemes. We show that there exist arithmetically Cohen-Macaulay
schemes that are not standard determinantal, and whose general
hyperplane section is good determinantal. We prove that if a general 
hyperplane section of a scheme is standard (resp. good)
determinantal, then the scheme is standard (resp. good) determinantal up to flat
deformation. We also study the transference of the property of being
standard or good determinantal under basic double linkage.
\end{abstract}

\section*{Introduction}

Standard and good determinantal schemes are a large family of
projective schemes, to which belong many varieties that have been
classically studied. For example rational
normal scrolls, rational normal curves, and some Segre varieties are
good determinantal schemes. Standard determinantal schemes are cut out
by the maximal minors of a matrix of polynomials (see
Definition~\ref{stdef}). In particular they are arithmetically
Cohen-Macaulay, and their saturated ideal is resolved by the
Eagon-Northcott complex. Good determinantal schemes are standard
determinantal schemes that are locally a complete
intersection outside a subscheme. 
Ideals of minors have been the object of extensive study
in commutative algebra. These families were studied from the geometric
viewpoint by Kreuzer, Migliore, Nagel, and Peterson in~\cite{kr00}. In
this article, they introduced the definition of standard and good determinantal
schemes that we use. The relevance of
standard and good determinantal schemes in the context of liaison
theory became clear in~\cite{kl01}, where it was
shown that standard determinantal schemes belong to the
Gorenstein-liaison class of a complete intersection. 

In this note we study standard and good determinantal schemes and
their general hyperplane sections. The property of being standard or
good determinantal is preserved when taking a general hyperplane
section. So we ask whether every arithmetically Cohen-Macaulay scheme
whose general hyperplane section is good determinantal is itself good
determinantal. The answer is negative. In Proposition~\ref{symm},
Example~\ref{vero}, and Proposition~\ref{n+1curve} we
produce examples of schemes which are not standard determinantal, and
whose general hyperplane section (or whose Artinian reduction) 
is good determinantal. We also show that a
section of the schemes of Proposition~\ref{symm} by a number of
generic hyperplanes is good determinantal up to flat deformation.
Then we discuss the property of being
standard or good determinantal in a flat family. This is motivated by
the observation that we can study flat families all of whose elements
are hyperplane section of a given scheme by a hyperplane that meets it 
properly. We show by means of
examples that we can have a flat family which contains a non standard
determinantal scheme and  whose general element is standard
determinantal, or the other way around. 
In Proposition~\ref{sect} we give sufficient conditions on a section of a
scheme $S$ by a hyperplane that meets it properly that force a
general hyperplane section of $S$ to be good determinantal. We saw that a
scheme $S$ with good determinantal general hyperplane section does not
need to be good determinantal. In Theorem~\ref{flatfam} we show that
$S$ is good determinantal up to flat deformation. Finally, 
we discuss how the property of being standard or good
determinantal is preserved under basic double linkage. In
Theorem~\ref{det} we prove that under some
assumptions the property is preserved. In Example~\ref{gensectbdl}
we show that in other cases the property is not preserved. We produce a
family of schemes via basic double link from the family of
Proposition~\ref{n+1curve}, and we prove that the schemes we produced
are not standard determinantal, but their general hyperplane sections
are good determinantal.

\section{Standard and good determinantal schemes}

Let $S$ be a scheme in $\PP^{n+1}=\PP^{n+1}_k$, where $k$ is an algebraically
closed field.
Let $I_S$ be the saturated homogeneous ideal corresponding to $S$ in the 
polynomial ring $R=k[x_0,\ldots,x_{n+1}]$. We denote by $\MM$ the
homogeneous irrelevant maximal ideal of $R$, $\MM=(x_0,\ldots ,x_{n+1})$.
For a sheaf $\F$ we denote by $$H^i_*(\F)=\bigoplus_{m\in\ZZ}H^i(\PP^{n+1},\F(m))$$ 
the i-th cohomology ring.
We will usually be interested in the case when $\F$ is an ideal sheaf.
Let $T$ be a scheme that contains $S$. We denote by $\I_{S|T}$ the
ideal sheaf of $S$ restricted to $T$, and by
$I_{S|T}=H^0_*(\I_{S|T})$ the ideal
of $S$ restricted to $T$.
We often write aCM for arithmetically Cohen-Macaulay.

In this note we study schemes whose general hyperplane section is
standard or good determinantal. The following definition was given
in~\cite{kr00} for schemes, i.e. for saturated ideals. 
Here we extend it to include Artinian ideals.
 
\begin{defn}\label{stdef}
An ideal $I\subseteq k[x_0,\ldots,x_{n+1}]$ of height $c$ 
is {\em standard determinantal} if it is 
generated by the maximal minors of a homogeneous matrix $M$
of polynomials of size $t\times (t+c-1)$, for some $t\geq 1$.
A matrix $M$ with polynomial entries is {\em homogeneous} if 
its minors are homogeneous polynomials.

A {\em standard determinantal} scheme 
$S\subseteq\PP^{n+1}$ of codimension $c$ is a scheme whose saturated
ideal $I_S$ is standard determinantal.

A standard determinantal ideal $I$ is {\em good determinantal} if 
after performing invertible row operations on the matrix $M$ and 
then deleting a row, the ideal generated by the maximal minors of 
the $(t-1)\times (t+c-1)$ matrix
obtained is standard determinantal (that is, it has height $c+1$). In
particular, we formally include the possibility that $t=1$, i.e. a
complete intersection is good determinantal.

A scheme $S$ is {\em good determinantal} if its saturated ideal $I_S$ is
good determinantal.
\end{defn}

Let $S$ be a standard determinantal scheme with defining matrix $M=(F_{ij})$. 
We assume without loss of generality that $M$ contains no invertible entries.
Let ${\mathcal U}=(u_{ji})$ be the transposed of the matrix whose entries 
are the degrees of the
entries of $M$. ${\mathcal U}$ is the {\it degree matrix} of $S$. We adopt the 
convention that the entries of ${\mathcal U}$ increase from right to left and 
from top
to bottom: $u_{ji} \geq u_{lk}$ if $i\leq k$ and $j \geq l$.
$S$ can be regarded as the degeneracy locus of a degree 0 
morphism $$\varphi:\bigoplus_{i=1}^t
R(b_i)\lra\bigoplus_{j=1}^{t+c-1} R(a_j).$$
Set $a_1\leq\ldots\leq a_{t+c-1}$ and $b_1\leq\ldots\leq b_t$. Then
$\varphi$ is described by the transposed of the matrix $M$, 
and $u_{ji}=a_j-b_i$.

\begin{defn}(\cite{ei88}, Definition~1.1)
A matrix $M=(F_{ij})$ is {\em 1-generic} if the entries in each row or
column are linearly independent over $k$.
\end{defn}

\begin{rmk}\label{1gen}
It is shown in~\cite{ei88} that the ideal generated by the maximal
minors of a 1-generic matrix
defines a reduced and irreducible standard determinantal scheme. Clearly
1-genericity is preserved if we delete a row of the
matrix. Therefore, the ideal of maximal minors of 
a 1-generic matrix defines a reduced and
irreducible good determinantal scheme.
\end{rmk}

In this note we study standard and good determinantal schemes, 
and schemes whose hyperplane section is standard or good determinantal.

\begin{defn}
Let $S\subseteq\PP^{n+1}$ be a projective scheme of dimension $d\geq 1$. 
Let $H$ be a hyperplane. If $H$ does not contain any component of $S$, we say that 
$S\cap H\subseteq H=\PP^n$ is a {\em proper hyperplane section} of $S$.

Fix a geometric property $\mathcal P$. We say that $\mathcal P$ holds
for a {\em general hyperplane section} of $S$ if there
is a nonempty open set $\mathcal V$ (in the $\PP^{n+1}$
parameterizing hyperplanes in $\PP^{n+1}$) such that $S\cap H$ has the
property $\mathcal P$ for all $H\in {\mathcal V}$. We call $S\cap
H\subseteq H=\PP^n$ a {\em general hyperplane section} of $S$.
\end{defn}

For a fixed scheme $S$, a general hyperplane section is proper. Namely, 
the set $\mathcal V$ of hyperplanes in $\PP^{n+1}$ that do not contain any 
component of $S$ is open and nonempty. 

If the scheme $S$ has dimension $d\geq 1$, then a general hyperplane section
has dimension $d-1$. 
If $I_S$ is the homogeneous saturated ideal of the scheme 
$S\subseteq\PP^{n+1}$, then $I_{S\cap H|H}=H^0_*(\PP^n,\I_{S\cap H})\subseteq R/(H)$ 
is the homogeneous saturated ideal of the
general hyperplane section $S\cap H\subseteq H$. 
The following short exact sequence of ideal sheaves
relates $S$ to a general hyperplane section $S\cap H$
$$0\lra\I_S(-1)\stackrel{\cdot H}{\lra}\I_S\lra\I_{S\cap H|H}\lra 0.$$
Taking cohomology we get the exact sequence:
$$\xymatrix{0\ar[r] & I_S(-1)\ar[r]^{\cdot H} & I_S\ar[rr]^{\pi}\ar[dr]& & I_{S\cap H|H}
\ar[r] & H^1_*(\I_S)(-1)\\
 & & & I_S+(H)/(H)\ar[ur]\ar[dr] & & \\
 & & 0\ar[ur] & & 0 & }$$
If $S$ is arithmetically 
Cohen-Macaulay, or more in general if $R/I_S$ has depth at least 2, 
then $H^1_*(\I_S)=0$, and $I_{S\cap H|H}=I_S+(H)/(H).$ 
In particular, if $S$ is arithmetically 
Cohen-Macaulay then a proper hyperplane section of $S$ is also arithmetically 
Cohen-Macaulay. Recall that a scheme $S$ of dimension $d\geq 0$ is
arithmetically Cohen-Macaulay if and only if $H^i_*(\I_S)=0$ for
$1\leq i\leq d$. Every zero-dimensional scheme is arithmetically Cohen-Macaulay.

If $S$ has dimension $d=0$, then geometrically it does not make sense to take 
a hyperplane section. However in this case the ideal
$I_S+(H)/(H)\subseteq R/(H)$ is 
Artinian (i.e. $R/I_S+(H)$ has Krull-dimension 0). In this case, 
we will abuse terminology and still call $I_S+(H)/(H)$ the ideal of a general 
hyperplane section of $S$, whenever $H\in R$ is a general linear form.
The short exact sequence relating the ideals of $S$ and of a general
hyperplane section is
$$0\lra I_S(-1)\stackrel{\cdot H}{\lra} I_S\lra I_S+(H)/(H)\lra 0.$$
We refer the interested reader to Section~1.3 of~\cite{mi98b} for
facts about hyperplane and hypersurface sections.

\section{Lifting the determinantal property, and good determinantal
  schemes in flat families}

In this note, we address the question of whether it is possible to lift the property 
of being standard or good determinantal from a general hyperplane
section of a scheme to the scheme itself. For schemes of codimension
$2$, the Hilbert-Burch Theorem states that being standard
determinantal is equivalent to being arithmetically Cohen-Macaulay. 
So this question is a natural generalization of the questions that were 
investigated by Huneke and Ulrich in~\cite{hu93}, by Migliore in~\cite{mi94}, 
and by the author in~\cite{go06}.

Before starting our discussion, we would like to observe that the
good determinantal property does not
behave as well as the standard determinantal property under hyperplane
sections by a hyperplane that meets the scheme properly. In fact,
any hyperplane section of a standard determinantal subscheme of
$\PP^{n+1}$ by a hyperplane that meets it properly is a
standard determinantal subscheme of $\PP^n$. It is not true in general
that every hyperplane section of a good determinantal
subscheme of $\PP^{n+1}$ by a hyperplane that meets it properly is a
good determinantal subscheme of $\PP^n$. However, a general
hyperplane section is good determinantal.
Next, we see an example when this is the case. The following example was derived 
from Example~4.1 in~\cite{kl01}.

\begin{ex}\label{stgood}
Let $C\subseteq\PP^4$ be a curve whose homogeneous saturated ideal is
given by the maximal minors of 
$$\left(\begin{array}{cccc}
x_0 & x_1+x_4 & 0 & x_2 \\
0 & x_1 & x_2 & x_0+x_1
\end{array}\right).$$
One can check that $C$ is one-dimensional, hence standard
determinantal. $C$ is a cone over a zero-dimensional scheme supported on 
the points $[0:0:0:1]$ and $[0:1:0:-1]$.
The curve $C$ is indeed good determinantal, since deleting a
generalized row we obtain the matrix of size $1\times 4$
$$\left(\begin{array}{cccc}
\alpha x_0 & (1+\alpha)x_1+\alpha x_4 & x_2 & x_0+x_1+\alpha x_2
\end{array}\right)$$
for a generic value of $\alpha$. For $\alpha\neq 0$ the entries form a
regular sequence, since they are 
linearly independent linear forms. Therefore they define a complete
intersection, that is a standard 
determinantal scheme, and $C$ is good determinantal.

Let $H$ be a general linear form. In particular we can assume that
the coefficient of $x_3$ in the equation of $H$ is
non-zero, so that $H$ does not contain the vertex of the cone $C$. 
Intersecting $C$ with $H$ we obtain a subscheme~$X$
of~$\PP^3$, whose saturated homogeneous ideal $I_X$ is generated over
$k[x_0,x_1,x_2,x_4]$ by the maximal minors of 
$$\left(\begin{array}{cccc}
x_0 & x_1+x_4 & 0 & x_2 \\
0 & x_1 & x_2 & x_0+x_1
\end{array}\right).$$
One can show that $X$ is good determinantal following the same steps
as for $C$. Indeed, $C$ is just a cone over $X$.

Let $H=x_4$. Intersecting $C$ with $H$ we obtain a subscheme $Z$ of $\PP^3$, 
whose saturated homogeneous ideal $I_Z$ is generated over
$k[x_0,\ldots,x_3]$ by the maximal minors of 
$$\left(\begin{array}{cccc}
x_0 & x_1 & 0 & x_2 \\
0 & x_1 & x_2 & x_0+x_1
\end{array}\right).$$
$I_Z=I_P^2$ for $P=[0:0:0:1]$, hence $Z$ is a zero-dimensional scheme
supported on the point $P$. Then $Z$ 
is standard determinantal and a section of $C$ by a hyperplane that
meets it properly. However, $Z$ is not good determinantal. In fact, deleting a
generalized row we obtain the matrix of size 
$1\times 4$
$$\left(\begin{array}{cccc}
\alpha x_0 & (1+\alpha)x_1 & x_2 & x_0+x_1+\alpha x_2
\end{array}\right)$$
whose entries generate the ideal $(x_0,x_1,x_2)$ of codimension
$3<4$.
\end{ex}
\vskip .5cm
Every standard determinantal scheme is arithmetically
Cohen-Macaulay. Moreover, the two families coincide for schemes of
codimension $1$ or $2$, while for codimension $3$ or higher the family
of arithmetically Cohen-Macaulay schemes strictly contains the family
of standard determinantal schemes.
From the results in \cite{hu93} one can easily obtain a sufficient condition
for a scheme $V\subseteq\PP^{n+1}$ to be arithmetically Cohen-Macaulay
in terms of the graded Betti numbers of a general hyperplane section
of $V$. If a general hyperplane section of $V$ is standard
determinantal, the condition can be expressed in terms of the entries
of its degree matrix. Notice that since the graded Betti numbers of a
hyperplane section of $V$ are the same for a general choice of the
hyperplane, the degree matrix is also the same for a general choice of the
hyperplane.

\begin{lemma}\label{aCM}
Let $V\subseteq\PP^{n+1}$ be a projective scheme. Assume that a
general hyperplane section of $V$ is a standard determinantal
subscheme of $\PP^n$ with degree
matrix\newline 
${\mathcal U}=(u_{ji})_{i=1,\ldots,t;\; j=1,\ldots,t+c-1}.$ If either 
$dim V\geq 2$ or
$$u_{1,t}+\cdots+u_{c-1,t}\geq n+1$$ then $V$ is arithmetically
Cohen-Macaulay.
\end{lemma}

\begin{proof}
If $dim(V)\geq 2$ and a general hyperplane section of $V$ is
arithmetically Cohen-Macaulay, then $V$ is arithmetically
Cohen-Macaulay (see Proposition~2.1
in~\cite{hu93}). We can then reduce to the
case when $V$ is one-dimensional.
Let $H$ be a general hyperplane, and let $Z=V\cap H$.
From Theorem~3.16 of~\cite{hu93} it follows that the minimum degree $b$ of
a minimal generator of $I_{Z|H}$ that is not the image of 
a minimal generator of $I_V$ under the standard projection
$I_V\stackrel{\pi}{\lra} I_{Z|H}$ is 
$$b\geq u_{1,1}+\cdots+u_{t,t}+u_{t+1,t}+\cdots+u_{t+c-1,t}-n=$$ 
$$=u_{1,t}+\cdots+u_{c-1,t}+u_{c,1}+u_{c+1,2}+\cdots+u_{t+c-1,t}-n
\geq u_{c,1}+u_{c+1,2}+\cdots+u_{t+c-1,t}+1.$$ 
In particular, it is bigger than the 
maximum $u_{c,1}+u_{c+1,2}+\cdots+u_{t+c-1,t}$ of the degrees of the
minimal generators of $I_{Z|H}$. Then 
all the minimal generators of $I_{Z|H}$ are images of the minimal
generators of $I_V$. Hence $H^1_*(\I_V)=0$, and 
$V$ is arithmetically Cohen-Macaulay.
\end{proof}

As we mentioned, every arithmetically
Cohen-Macaulay scheme of codimension $2$ is standard determinantal. So
Lemma~\ref{aCM} gives a sufficient condition to conclude that $V$
is standard determinantal if $codim(V)=2$.

\begin{rmk}\label{dim2}
Let $V$ be a projective scheme. If $dim(V)\geq 2$ and a general
hyperplane section of $V$ is aCM, then $V$ is aCM. 
Therefore the graded Betti numbers of $V$ coincide with the graded Betti 
numbers of a general hyperplane section of $V$ (for more details 
see~\cite{mi98b}, Theorem~1.3.6). Moreover,
for a scheme of codimension $2$ the property of being standard
determinantal can be decided by checking the graded Betti
numbers. In fact, a scheme of codimension $2$ is standard determinantal if and only
if it is aCM, if and only if a minimal free
resolution of its saturated ideal has length $2$. 
Hence if $dim(V)\geq 2$ and $codim(V)=2$, we can decide
whether $V$ is standard determinantal by looking at the graded Betti
numbers of a general hyperplane section.
However, if $codim(V)\geq 3$ then
the property of being standard determinantal cannot in general be
decided by looking at the graded Betti numbers. In other words, there
are schemes which are not standard determinantal, but have the same
graded Betti numbers as a standard determinantal scheme 
(see e.g. Example~\ref{vero}).
\end{rmk}

In very special cases the graded Betti numbers of a homogeneous ideal
$I$ can force the ideal to be standard determinantal, even when the
codimension is $3$ or higher. The next is an easy example of this phenomenon.

\begin{ex}\label{artin}
Let $R=k[x_1,\ldots,x_n]$, $\MM=(x_1,\ldots,x_n)$. Let
$I\subseteq R$ be a homogeneous ideal generated by 
$\bin{n+t-1}{t}$ linearly independent polynomials of degree $t$.
Then $I_j=0$ for all $j<t$ and
$dim I_t=\bin{n+t-1}{t}=dim(\MM^t)_t$. 
Therefore 
$I=\MM^t$, so it is the ideal of 
maximal minors of the $t\times (t+n-1)$ matrix $$\left(\begin{array}{ccccccc}
x_1 & \cdots & x_n & 0 & \cdots & \cdots & 0 \\
0 & x_1 & \cdots & x_n & 0 & & \vdots \\
\vdots & \ddots & \ddots & \ddots & \ddots & \ddots & \vdots \\
\vdots &  & 0 & x_1 & \cdots & x_n & 0 \\
0 & \cdots & \cdots & 0 & x_1 & \cdots & x_n 
\end{array}\right).$$
So $I$ is good determinantal. 
\end{ex}

The next proposition shows that this is not the case in general. We
present a family of arithmetically Cohen-Macaulay schemes that are not standard
determinantal, but such that the Artinian reduction of their
coordinate ring is good determinantal.
In particular, they have the graded Betti numbers of a
standard determinantal scheme. 
From a more geometric point of view, it is interesting to decide whether
the schemes in question have a general section which is 
good determinantal. In other words, whether a
section of $V$ by $r$ generic hyperplanes is good
determinantal for some $r\leq {t+2\choose 2}-3$. We prove that the
schemes of the following proposition have a (special) $\bin{t}{2}$-th
proper hyperplane section
which is good determinantal.

\begin{prop}\label{symm}
Let $X$ be a symmetric matrix of indeterminates of size
$(t+1)\times(t+1)$, $t\geq 2$ 
$$X=\left(\begin{array}{ccccc} x_{0,0} & x_{0,1} & \cdots & \cdots & x_{0,t} \\
x_{0,1} & x_{1,1} & \cdots & \cdots & x_{1,t} \\
\vdots & \vdots & & & \vdots \\
x_{0,t} & x_{1,t} & \cdots & \cdots & x_{t,t} \end{array}\right).$$
Let $V\subseteq\PP^{\bin{t+2}{2}-1}$ be the scheme corresponding
to the saturated ideal 
$I_V=I_t(X)\subseteq R=k[\; x_{i,j}\; |\; 0\leq i\leq j\leq t\;]$,
generated by the submaximal 
minors of $X$.
\begin{enumerate}
\item $V$ is an arithmetically Cohen-Macaulay, integral scheme of codimension $3$ which 
  is not standard determinantal, but every Artinian reduction of its
  homogeneous coordinate ring is good determinantal.
\item Let $D$ be a general $\bin{t}{2}$-th hyperplane section of
  $V$. Then $V$ has a proper $\bin{t}{2}$-th hyperplane section $C$ that is a
  good determinantal scheme, and there is a flat family of schemes
  with fixed graded Betti numbers that contains both $C$ and $D$.
\end{enumerate}
\end{prop}

\begin{proof}
(1) The fact that $V$ is an arithmetically Cohen-Macaulay, integral scheme of 
codimension $3$ follows from classical results, that can be found e.g. 
in~\cite{br88}. 
In particular a minimal free resolution of the ideal $I_V$ is known, and
the cardinality of a minimal system of generators of
$I_V$ is $m=\bin{t+2}{2}$. Then the Artinian reduction of the coordinate 
ring of $V$ is good determinantal, as showed in Example~\ref{artin}. 
The divisor class group of $V$
is isomorphic to $\ZZ_2$ (see~\cite{go77}).
From knowledge of the graded Betti numbers of $I_V$ (see
e.g.~\cite{br88}), it follows that
if $V$ was standard determinantal, then its degree matrix would have
size $t\times (t+2)$ and all of its entries would be equal to $1$. 
The divisor class group of such a standard determinantal scheme is
isomorphic to $\ZZ$ (see~\cite{br75}). Therefore $V$ is not standard
determinantal. 
Notice that if $t=2$ then $V$ is the Veronese surface in
$\PP^5$, which is not standard determinantal, since it is not
isomorphic to a rational normal scroll surface. 

(2) Consider a special $\bin{t}{2}$-th hyperplane section of $V$,
with defining matrix of size $(t+1)\times (t+1)$
$$Y=\left(\begin{array}{cccccc} 
x_{0,0} & x_{0,1} & x_{0,2} & \cdots & \cdots & x_{0,t} \\ 
x_{0,1} & x_{0,2} & & & x_{0,t} & x_{1,t} \\
x_{0,2} & & & x_{0,t} & x_{1,t} & \vdots \\
\vdots & & \adots & \adots & & \vdots \\
\vdots & x_{0,t} & x_{1,t} & & & x_{t-1,t} \\
x_{0,t} & x_{1,t} & \cdots & \cdots & x_{t-1,t} & x_{t,t} 
\end{array}\right).$$
We obtain this section intersecting with the hyperplanes 
$x_{i,j}-x_{0,i+j}$ for $i+j\leq t$ and $i\geq 1, j\leq t-1$ and
$x_{i,j}-x_{i+j-t,t}$ for $i+j>t$ and $i\geq 1, j\leq t-1$.
We take $\bin{t}{2}$ hyperplane sections by hyperplanes that meet $V$
properly. So we obtain a scheme
$C\subseteq\PP^{2t}$ of codimension $3$.
$C$ is good determinantal, with defining matrix 
$$U=\left(\begin{array}{cccccccc} 
x_{0,0} & x_{0,1} & x_{0,2} & \cdots & \cdots & x_{0,t-1} & x_{0,t} & x_{1,t} \\ 
x_{0,1} & x_{0,2} & & & x_{0,t-1} & x_{0,t} & x_{1,t} & x_{2,t} \\
x_{0,2} & & & x_{0,t-1} & x_{0,t} & x_{1,t} & x_{2,t} & \vdots \\
\vdots & & \adots & \adots & \adots & \adots & & \vdots \\
\vdots & x_{0,t-1} & x_{0,t} & x_{1,t} & x_{2,t} & & & \vdots \\ 
x_{0,t-1} & x_{0,t} & x_{1,t} & x_{2,t} & \cdots & \cdots & \cdots & x_{t,t} 
\end{array}\right).$$
In fact, the maximal minors of $U$ coincide with the submaximal minors
of $Y$. Moreover, the matrix $U$ is 1-generic.  
Therefore, the ideal of maximal minors of $U$ defines a reduced
and irreducible, good determinantal scheme (see also Remark~\ref{1gen}). 

Let $D$ be a general $\bin{t}{2}$-th hyperplane section of $V$. 
The saturated ideal of $D$ is the ideal $I_D=I_t(Z)$ generated
by the submaximal minors of the symmetric matrix 
$$Z=\left(\begin{array}{ccccc} 
x_{0,0} & x_{0,1} & \cdots & x_{0,t-1} & x_{0,t} \\ 
x_{0,1} & L_{1,1} & \cdots & L_{1,t-1} & x_{1,t} \\
\vdots & \vdots & & \vdots & \vdots \\
x_{0,t-1} & L_{1,t-1} & \cdots & L_{t-1,t-1} & x_{t-1,t} \\
x_{0,t} & x_{1,t} & \cdots & x_{t-1,t} & x_{t,t} \end{array}\right).$$
We can assume without loss of generality that the equations of the 
hyperplanes that we intersect with $V$ are $x_{i,j}-L_{i,j}$, 
$i\geq 1, j\leq t-1$, where
$L_{i,j}$ is a general linear form in $k[x_{0,0},\ldots, x_{0,t},
x_{1,t},\ldots, x_{t,t}]$. Observe that we
have a flat family of codimension 3 schemes $D_s$ whose saturated ideal is
$I_t(Z_s)$, $Z_s=sZ+(1-s)Y$. In fact, for any choice of $s$ and for
$L_{i,j}$ generic, the matrix $Z_s$ is 1-generic. 
Then by Corollary~3.3 of \cite{ei88}  
$$codim\: I_t(Z_s)\geq 2(t+1)-1-2(t-1)=3.$$
Hence $Z_s$ defines an aCM scheme $D_s$ of codimension three, whose
graded Betti numbers are the same as those of $C$ and of $V$ (this
follows from~\cite{br88}, Theorem~3.5). In
particular the Hilbert polynomial of $D_s$ is the same
for all $s$. 
\end{proof}

The Veronese surface $V\subseteq\PP^5$ is an example of a non standard
determinantal scheme from the family of
Proposition~\ref{symm}.
In the next example we show that a general hyperplane section of $V$ is
a good determinantal curve.

\begin{ex}\label{vero}
The Veronese surface $V\subseteq\PP^5$ is an example from the family of
Proposition~\ref{symm}, for $t=2$.
Its homogeneous saturated ideal is the ideal 
$$I_V=I_2\left(\begin{array}{ccc} x_0 & x_1 & x_2 \\
x_1 & x_3 & x_4 \\ x_2 & x_4 & x_5 \end{array}\right)$$ 
$I_V\subseteq S=k[x_0,\ldots,x_5].$
Its general hyperplane section is a reduced and irreducible
arithmetically Cohen-Macaulay curve 
$C\subseteq\PP^4$ of degree $4$, hence a rational normal
curve. In particular, a general hyperplane section of $V$ is good
determinantal, with defining matrix equal to (after a
change of coordinates and invertible row and column operations)
$$\left(\begin{array}{cccc}  x_0 & x_1 & x_2 & x_3 \\
x_1 & x_2 & x_3 & x_4 \end{array}\right).$$
\end{ex}
\vskip .5cm
Kleppe, Migliore, Mir\'o-Roig, Nagel and
Peterson proved that under certain assumptions the closure of the locus of good
determinantal schemes with a fixed degree matrix $M$ is an irreducible
component in the corresponding Hilbert scheme (see chapters~9 and~10
of~\cite{kl01} and the paper \cite{kl05}). 
Clearly, standard determinantal schemes with the same degree matrix
$\mathcal U$ belong to the closure of the locus of good determinantal
ones. It is natural to ask whether a general ${t\choose 2}$-th
hyperplane section of a scheme $V$
as in Proposition~\ref{symm} is standard (or good) determinantal. The following
example shows that this is in general not the case.

\begin{ex}\label{verodeform}
Let $V\subseteq\PP^9$ be the scheme whose saturated homogeneous ideal
$I_V$ is generated by the submaximal minors of
the matrix $$X=\left(
\begin{array}{cccc}
x_{0,0} & x_{0,1} & x_{0,2} & x_{0,3} \\
x_{0,1} & x_{1,1} & x_{1,2} & x_{1,3} \\
x_{0,2} & x_{1,2} & x_{2,2} & x_{2,3} \\
x_{0,3} & x_{1,3} & x_{2,3} & x_{3,3}
\end{array}\right).$$
In Proposition~\ref{symm} we showed that $V$ has a $3$-rd hyperplane
section $C\subseteq\PP^6$ that is good determinantal. More precisely, the
ideal $I_C$ is generated by the maximal minors of the matrix 
$$Y=\left(\begin{array}{ccccc}
x_{0,0} & x_{0,1} & x_{0,2} & x_{0,3} & x_{1,3} \\
x_{0,1} & x_{0,2} & x_{0,3} & x_{1,3} & x_{2,3} \\
x_{0,2} & x_{0,3} & x_{1,3} & x_{2,3} & x_{3,3}
\end{array}\right).$$
The homogeneous saturated ideal of a general
$3$-rd hyperplane section of $V$ is generated by the maximal minors of
the matrix 
$$Z=\left(
\begin{array}{cccc}
x_{0,0} & x_{0,1} & x_{0,2} & x_{0,3} \\
x_{0,1} & L_{1,1} & L_{1,2} & x_{1,3} \\
x_{0,2} & L_{1,2} & L_{2,2} & x_{2,3} \\
x_{0,3} & x_{1,3} & x_{2,3} & x_{3,3}
\end{array}\right)$$
where $L_{1,1},L_{1,2},L_{2,2}\in
k[x_{0,0},x_{0,1},x_{0,2},x_{0,3},x_{1,3},x_{2,3},x_{3,3}]$ 
are general linear forms.
Let $I(s)=I_3(Z_s)$ be the ideal generated by the submaximal minors of
the matrix
$Z_s=sZ+(1-s)Y$. Then one can check
that for a generic value of $s$
the cardinality of a minimal system of generators of $I(s)^2$ is 
$\mu(I(s)^2)=55$ (we used the computer algebra software CoCoA
\cite{cocoa}). 
If $I(s)$ defines a standard determinantal scheme, then it follows by
knowledge of the graded Betti numbers of $I(s)$ that it must be
associated to a matrix of linear forms of size $3\times 5$.  
In that case we have 5 linearly independent Pl\"ucker relations, which
implies that $\mu(I(s)^2)\leq 50$. Therefore $I(s)$ cannot define a
standard determinantal scheme. 
Hence $V$ has a good determinantal $3$-rd hyperplane section by
hyperplanes that meet it properly, while its general $3$-rd
hyperplane section is not standard determinantal.
\end{ex} 

The last family of examples that we wish to study consists of non 
standard determinantal curves, whose general hyperplane section 
is good determinantal (see Proposition~\ref{n+1curve}). 
The result of the next lemma is not new. 
For completeness we give a simple algebraic proof of it.

\begin{lemma}\label{sqfr}
Let $I\subseteq k[x_0,\ldots,x_n]$ be the ideal generated by all the
squarefree monomials of degree $d$. Then $I$ is a good determinantal 
ideal.
\end{lemma}

\begin{proof}
Let $A$ be a matrix of size $d\times (n+1)$ with entries in $k$ such that
all the maximal minors of $A$ are nonzero, $A=(a_{i,j})_{1\leq i\leq
  d;\; 0\leq j\leq n}$.
Consider the matrix $M$ that we obtain from $A$ by multiplying each
entry in the $j$-th column by $x_j$, $M=(a_{i,j}x_j)_{1\leq i\leq
  d;\; 0\leq j\leq n}.$ 
The minor involving columns $0\leq j_1<\ldots<j_d\leq n$ is
$\alpha_{j_1,\ldots,j_d} x_{j_1}\cdot\ldots\cdot x_{j_d}$, where
$\alpha_{j_1,\ldots,j_d}$ is the determinant of the submatrix of $A$
consisting of the columns $j_1,\ldots,j_d$.
If $d=1$ then $I$ is a complete intersection, hence good determinantal.
If $d\geq 2$ the height of $I$ is $n+2-d$, 
then $I$ is standard determinantal. In particular,
$k[x_0,\ldots,x_n]/I$ is Cohen-Macaulay. If
we delete a generalized row of $M$, up to nonzero scalar multiples 
the $(d-1)\times (d-1)$ minors of
the remaining rows are all the squarefree monomials of degree $d-1$ in 
$k[x_0,\ldots,x_n]$. Since they generate a standard determinantal
ideal, $I$ is good determinantal.
\end{proof}

\begin{rmk}
In order for the result of Lemma~\ref{sqfr} to hold we do not even need the
ground field $k$ to have infinite cardinality. 
However we
need to have enough scalars in $k$ so that we can find a matrix $A$ 
of size $d\times (n+1)$ with entries in $k$ such that
all the $d\times d$ minors of $A$ are nonzero. If $|k|\geq
d+1$ we can let $A$ be the Vandermonde matrix in $\alpha_1,\ldots,\alpha_d$,
distinct elements in $k^*$, i.e. $a_{ij}=\alpha_i^{j-1}$.
\end{rmk}

The next proposition is a straightforward consequence of Lemma~\ref{sqfr}.

\begin{prop}\label{genpts}
$n+1$ generic points in $\PP^n$ are a good determinantal scheme.
\end{prop}

\begin{proof}
Observe that $n+1$ generic points in $\PP^n$ can be mapped via
a change of coordinates to the $n+1$ coordinate points. 
The saturated ideal of the $n+1$ coordinate points in $\PP^n$ is
generated by the squarefree monomials of degree 2 in $x_0,\ldots,x_n$.
Therefore it is a good determinantal scheme by Lemma~\ref{sqfr}. 
\end{proof}

Let $C\subseteq\PP^{n+1}$ be a nondegenerate, reduced and irreducible
curve of degree $n+1$. Then $C$ is a rational normal curve, in
particular it is good determinantal.
In the next proposition we produce a nondegenerate, 
arithmetically Cohen-Macaulay, 
reduced curve of degree $n+1$ in
$\PP^{n+1}$ that is not standard determinantal and whose general
hyperplane section is good determinantal. The curve is necessarily
reducible, because of what we just observed.

\begin{prop}\label{n+1curve}
Let $C_1\subseteq\PP^{n+1}$ be a cone over $n$ generic points in
$\PP^n$. Let $C_2\subseteq\PP^{n+1}$ be a generic line through a point
in $C_1$. Let $C=C_1\cup C_2$.
Then $C$ is not standard determinantal, and a general hyperplane
section of $C$ is good determinantal.
\end{prop}

\begin{proof}
Without loss of generality we can let the $n$ generic points in
$\PP^n$ be all the coordinate points except for $[1:0:\ldots:0]$.
Then the saturated ideal of $C_1\subseteq\PP^{n+1}$ is 
$$I_{C_1}=\bigcap_{i=1}^n(x_0,x_1,\ldots,x_{i-1},x_{i+1},\ldots,x_n)=
(x_0)+\sum_{1\leq i<j\leq n}(x_ix_j).$$
We can also assume that
$I_{C_2}=(x_2,\ldots,x_{n+1}).$
Then the saturated ideal of $C$ is
$$I_C=I_{C_1}\cap I_{C_2}=x_0(x_2,\ldots,x_{n+1})
+\sum_{1\leq i<j\leq n}(x_ix_j).$$
Since $I_{C_1}+I_{C_2}=(x_0,x_2,\ldots,x_{n+1})=I_P$ where $P$ is the
point $[0:1:0:\ldots:0]$, it follows that $C$ is arithmetically
Cohen-Macaulay. This follows from computing cohomology from
the short exact sequence
$$0\lra\I_C\lra \I_{C_1}\oplus\I_{C_2}\lra \I_P\lra 0.$$
In fact we obtain the long exact sequence
$$0\lra I_C\lra I_{C_1}\oplus I_{C_2}\lra
I_{C_1}+I_{C_2}=I_P\stackrel{0}{\lra} H^1_*(\I_C)\lra
H^1_*(\I_{C_1})\oplus H^1_*(\I_{C_2})=0$$
where vanishing of the last module follows from the observation that
$C_1$ and $C_2$ are aCM.
The curve $C$ has degree $n+1$, and its general
hyperplane section consists of $n+1$ generic points in $\PP^n$ by
construction. 
Therefore a general
hyperplane section of $C$ is good determinantal by Proposition~\ref{genpts}.

We now study the last morphism in a minimal free resolution of $I_C$,
in order to show that $C$ is not standard determinantal.
Let $I=x_0(x_2,\ldots,x_{n+1})$, $J=\sum_{1\leq i<j\leq
  n}(x_ix_j)$. Then clearly $I_C=I+J$. So we have the short exact
sequence
\begin{equation}\label{mapp}
0\lra I\cap J\lra I\oplus J\lra I+J\lra 0.
\end{equation}
Let \begin{equation}\label{mfrI}
0\lra\FF_n\lra\FF_{n-1}\lra\cdots\lra\FF_1\lra I\lra 0
\end{equation}
be a minimal free resolution of $I$. 
Then $\FF_n=R(-n-1)$ and $\FF_{n-1}=R(-n)^n.$
The last morphism in (\ref{mfrI}) is 
\begin{equation}\label{lastI}
R(-n-1)\stackrel{(x_2,-x_3,x_4,\ldots,(-1)^{n+1}x_{n+1})}{\lra}
R(-n)^n.\end{equation}
Let 
\begin{equation}\label{mfrJ}
0\lra\GG_{n-1}\stackrel{M}{\lra}\GG_{n-2}\lra\cdots\lra\GG_1\lra J\lra 0
\end{equation}
be a minimal free resolution of $J$. 
The ideal $J$ is a lexsegment squarefree monomial ideal, hence morphisms in a
minimal free resolution are explicitly computed in \cite{ar98},
Theorem~2.1. It turns out that $\GG_{n-1}=R(-n)^{n-1}$, 
$\GG_{n-2}=R(-n+1)^{n(n-2)}$, and
the matrix $M$ describing the last morphism in (\ref{mfrJ}) has size 
$n(n-2)\times(n-1)$ and is of the form 
$$M=\left(\begin{array}{ccc}c_1 & \ldots & c_{n-1} 
\end{array}\right)$$ where each $c_i$ is a column with exactly
$n-1$ nonzero entries (all the indeterminates but $x_i$). 
Finally, $I\cap J=x_0J$, so the minimal free resolution (\ref{mfrJ})
twisted by $-1$ is a minimal free resolution of $I\cap J$
\begin{equation}\label{mfrIJ}
0\lra\GG_{n-1}(-1)\stackrel{M}{\lra}\GG_{n-2}(-1)\lra\cdots\lra\GG_1(-1)\lra I\cap J\lra 0.
\end{equation}

Using the mapping cone construction on the short exact sequence
(\ref{mapp}), one can write the last matrix in a minimal free
resolution of $I+J=I_C$:
\begin{equation}\label{matx}\left(\begin{array}{ccccc}
x_2 & 0 & \ldots & \ldots & 0 \\
-x_3 & -x_1 & 0 & \ldots & 0 \\
\vdots & 0 & \ddots & & \vdots \\
\vdots & \vdots & & \ddots & 0 \\
(-1)^{n+1} x_{n+1} & 0 & \ldots & 0 & -x_1 \\
0 & x_0 &  0 & \ldots & 0 \\
\vdots & 0 & \ddots & & \vdots \\
\vdots & \vdots & & \ddots & 0 \\
0 & 0 & \ldots & 0 & x_0 \\
0 & & & & \\
0 & & M & & \\
0 & & & & \\
\end{array}\right).\end{equation}
The matrix corresponds to a morphism 
$$R(-n-1)^n=\FF_n\oplus\GG_{n-1}(-1)\lra 
\FF_{n-1}\oplus\GG_{n-1}\oplus\GG_{n-2}(-1)=R(-n)^{n^2-1},$$
where the block consisting of the first $n-1$ rows and the first
column comes from the last map in a minimal free resolution of $I$,
i.e. (\ref{lastI}).
The block consisting of the last $n(n-2)$ rows and
the last $n-1$ columns comes from the last map $M$ in a minimal free
resolution of $I\cap J$. The block consisting of the first $2n-1$ rows and
last $n-1$ columns comes from the morphism 
$\GG_{n-1}(-1)\lra\FF_{n-1}\oplus\GG_{n-1}$ induced by the diagonal morphism
$I\cap J\lra I\oplus J$.
The rows $2,\ldots,n$ have $-x_1$ on the diagonal and zeroes
anywhere else, while the rows $n+1,\ldots,2n-1$ have 
$x_0$ on the diagonal and zeroes
anywhere else. This corresponds to the fact that each minimal
generator of $I\cap J$ is of the form $x_0$ multiplied by a minimal
generator of $J$, which is also equal to $x_1$ multiplied by a minimal
generator of $I$. The indeterminate $x_{n+1}$ appears in the matrix 
(\ref{matx}) only in one position. From this observation and from the
form of $M$ it is easy to see that the ideal of $2\times 2$ minors of the matrix 
(\ref{matx}) is $(x_0,\ldots,x_n)^2+x_{n+1}(x_0,\ldots,x_n).$

Suppose by contradiction that the curve $C$ is standard determinantal. Then there
exist linear forms $L_0,\ldots,L_{2n+1}\in k[x_0,\ldots,x_{n+1}]$ such
that $I_C$ is the ideal of $2\times 2$ minors of the matrix 
$$\left(\begin{array}{ccc}
L_0 & \ldots & L_n \\
L_{n+1} & \ldots & L_{2n+1}
\end{array}\right).$$
The Eagon-Northcott complex is a minimal free resolution of the ideal
$I_C$. The last matrix in the complex has a block form, where
the basic block is given by the two column vectors
$$U=\left(\begin{array}{c}
-L_{n+1} \\
L_{n+2} \\
\vdots \\
(-1)^{n+1} L_{2n+1}
\end{array}\right)\;\;\;\mbox{and}\;\;\; V=\left(\begin{array}{c}
L_0 \\
-L_1 \\
\vdots \\
(-1)^n L_n
\end{array}\right).$$ 
The matrix has the form 
\begin{equation}\label{matx2}
\left(\begin{array}{cccccc}
U & V & 0 & 0 & \ldots & 0 \\
0 & U & V & 0 & \ldots & 0 \\
\vdots & \ddots & \ddots & \ddots & \ddots & \vdots \\
0 & \ldots & 0 & U & V & 0 \\ 
0 & \ldots & 0 & 0 & U & V 
\end{array}\right).
\end{equation}
In particular the ideal of $2\times 2$ minors of the matrix 
(\ref{matx2}) is
$$(L_0,\ldots,L_{2n+1})^2=(x_0,\ldots,x_n)^2+x_{n+1}(x_0,\ldots,x_n).$$
Taking radicals we obtain
$$(L_0,\ldots,L_{2n+1})
\subseteq\sqrt{(L_0,\ldots,L_{2n+1})^2}=(x_0,\ldots,x_n)$$ hence
$x_{n+1}\not\in (L_0,\ldots,L_{2n+1})$, which is a contradiction.
\end{proof}

In Proposition~\ref{symm}, Example~\ref{vero}, and
Proposition~\ref{n+1curve}
we discussed some examples of ``pathological'' behavior connected with lifting
the property of being standard or good determinantal. The schemes we
studied are all defined by minors of matrices with linear entries.
In analogy with the question of lifting the property of being
arithmetically Cohen-Macaulay (see Lemma~\ref{aCM}), 
one could ask the following.

\begin{quest}\label{asympt}
Assume that $char(k)=0$ and let $V\subseteq\PP^{n+1}_k$ be an aCM
scheme. Let $C\subseteq\PP^n$ be a general hyperplane section of $V$,
and assume that $C$ is standard/good determinantal. Does there exist
an $N$ such that if all the entries of the degree matrix of $C$ are at
least $N$, then $V$ is standard/good determinantal?
\end{quest}

The next example illustrates the necessity of requiring that a general
hyperplane section of the scheme is standard (or good) determinantal,
as opposed to requiring that a hyperplane section by a hyperplane that
meets the scheme properly is standard (or good) determinantal. 
Notice that the entries of the degree 
matrix $M$ in the next example can be taken arbitrarily large.

\begin{ex}
Let $k$ have arbitrary characteristic.
Let $V\subseteq\PP^5_k$ be the scheme corresponding to the saturated ideal
$$I_V=I_2\left(\begin{array}{ccc}
x_0^n & x_1^n & x_2^n \\
x_1^n & x_3^n & x_4^n \\
x_2^n & x_4^n & x_5^n
\end{array}\right)\subseteq k[x_0,\ldots,x_5].$$
The ideal $I_V$ is saturated and has height $3$, hence it defines a surface 
$V\subseteq\PP^5$. Since $\hgt I_V=3$, a minimal free resolution of $I_V$ 
can be obtained from a minimal free resolution of the Veronese surface 
by substituting $x_i$ by $x_i^n$ for $i=0,\ldots,5$. This follows 
from Theorem~3.5 in~\cite{br88}. One can check that $V$ is not standard 
determinantal by a similar argument to that used for the Veronese surface 
in Proposition~\ref{vero}. 

Let us intersect $V$ with a hyperplane $H$ of equation
$x_3-x_2=0$. 
The scheme $D=V\cap H$ is arithmetically Cohen-Macaulay, and its 
saturated ideal $I_D$ is generated by the submaximal minors of the matrix
$$\left(\begin{array}{ccc}
x_0^n & x_1^n & x_2^n \\
x_1^n & x_2^n & x_4^n \\
x_2^n & x_4^n & x_5^n
\end{array}\right).$$
Consider a rational normal curve $C$ whose saturated ideal $I_C$ 
is generated by the submaximal minors of the matrix 
$$\left(\begin{array}{ccc}
x_0 & x_1 & x_2 \\
x_1 & x_2 & x_4 \\
x_2 & x_4 & x_5
\end{array}\right).$$
$C$ is good determinantal, hence the Eagon-Northcott complex is a minimal 
free resolution of $I_C$. 
Since $\hgt I_D=\hgt I_C=3$ and $I_D$ is obtained from $I_C$ by replacing each 
occurrence of $x_i$ by $x_i^n$, it follows from Theorem~3.5 in~\cite{br88} 
that we can obtain a minimal free resolution of $I_D$ from a minimal free 
resolution of $I_C$ by replacing each occurrence of $x_i$ by $x_i^n$.
$D$ is good determinantal, since $C$ is.
\end{ex}

We now present an easy example that shows how the closure of the locus of
good determinantal schemes in the Hilbert scheme can contain also
schemes that are not standard determinantal (or not even
arithmetically Cohen-Macaulay).

\begin{ex}\label{deg9gen10}
Consider the Hilbert scheme $\HH$ parameterizing curves of degree $9$
and genus $10$ in $\PP^3$.
Let $D$ be the locus of $\HH$ whose points correspond to a
$CI(3,3)$. Let $E$ be the locus of $\HH$ whose points correspond to
curves of type $(3,6)$ on a smooth quadric surface.
The elements of $E$ are non-aCM. In fact, up to linear equivalence, a
curve of type $(3,6)$ is $C=C_1\cup C_2$ where $C_1$ consists of $3$
skew lines and $C_2$ consists of $6$ skew lines. Moreover, each line
of $C_1$ intersects each line of $C_2$, so $C_1\cap C_2$ consists of
$18$ distinct points. Let $I_C\subset R=k[x_0,x_1,x_2,x_3]$ be the
ideal corresponding to $C$. The minimal free resolution of $I_C$ as
an $R$-module is 
$$0\lra R^2(-8)\lra R^6(-7)\lra R^4(-6)\oplus R(-2)\lra I_C\lra 0.$$
In particular, $C$ is non-aCM.

By the uppersemicontinuity principle, no point of the closure of $E$
can be aCM, so $E$ is closed. But since $\HH$ is connected, the
closure of $D$ needs to intersect $E$, therefore there is a point in
the closure of $D$ that corresponds to a non-aCM curve. Notice that
aCM schemes and standard determinantal schemes coincide in the
codimension $2$ case. So this shows that the closure of the locus of
good determinantal schemes in the Hilbert scheme can contain also
schemes that are not standard determinantal (and not even
arithmetically Cohen-Macaulay). 
\end{ex}

Examples~\ref{verodeform} and~\ref{deg9gen10} show that we can have a
flat family which contains a non standard determinantal scheme and 
whose general element is standard determinantal, or the other way
around. Notice however that while all the schemes in the flat family
of Example~\ref{verodeform} are arithmetically Cohen-Macaulay, the non
standard determinantal element in the flat family of
Example~\ref{deg9gen10} is not aCM.

In Example~\ref{verodeform} we exhibit an arithmetically
Cohen-Macaulay scheme that has a proper 3-rd hyperplane section which
is good determinantal, but whose general 3-rd hyperplane section is
not good determinantal. Under some assumptions we can conclude that 
if a scheme $V$ has a good determinantal section by a hyperplane that
meets $V$ properly, then a general hyperplane section of $V$ is good 
determinantal. In the sequel, we will see that this forces $V$ to be 
good determinantal up to flat deformation (see Theorem~\ref{flatfam}).

Let $S\subseteq\PP^{n+1}$ be a scheme of dimension $d\geq 2$ and let
$C\subseteq\PP^n$ be a general hyperplane section of $S$.
Notice that since we are working with schemes of dimension greater
than or equal to 1, it is not restrictive to
assume that $S$ is arithmetically Cohen-Macaulay. In fact, $C$ aCM
forces $S$ to be aCM. Sufficient conditions for the unobstructedness
of $C$ are discussed in the last two chapters of~\cite{kl01}.

\begin{prop}\label{sect}
Let $k$ have characteristic zero.
Let $S\subseteq\PP^{n+1}$ be an aCM scheme and let $C\subseteq\PP^n$
be a hyperplane section of $S$ by a hyperplane that meets $S$
properly. Assume that $C$ is good determinantal, and let ${\mathcal U}=(u_{ij})$ be the
degree matrix of $C$. Let $p$ be the Hilbert polynomial of $S$ and
$C$, and let $Hilb^p(\PP^n)$ be the Hilbert scheme of subschemes of
$\PP^n$ with Hilbert polynomial $p$. Assume that $C$ belongs to the interior of the
locus of good determinantal schemes with degree matrix ${\mathcal U}$ 
in $Hilb^p(\PP^n)$, and that it is unobstructed. 
Assume moreover that one of the following holds:
\begin{itemize}
\item $S$ has codimension $3$, and $n\geq 5$;
\item $S$ has codimension $3$, $n\geq 4$, $u_{i,i-\mmin\{2,t\}}\geq 0$
  for $\mmin\{2,t\}\leq i\leq t$, and $u_{t,t+1}>u_{t,t}+u_{1,t-1}$;
\item $S$ has codimension $3$, $n=4$, and $u_{t,0}>u_{t,1}+u_{t,2}$;
\item $S$ has codimension $4$, $n\geq 6$, and $u_{i,i-\mmin\{3,t\}}\geq 0$
  for $\mmin\{3,t\}\leq i\leq t$;
\item $S$ has codimension $4$, $n\geq 5$, $u_{i,i-\mmin\{3,t\}}\geq 0$
  for $\mmin\{3,t\}\leq i\leq t$, and $u_{t,t+2}>u_{t,t}+u_{1,t-1}$;
\item $S$ has codimension $c\geq 5$, $n\geq c+1$, $u_{i,i-\mmin\{3,t\}}\geq 0$
  for $\mmin\{3,t\}\leq i\leq t$, and $u_{t,t+j-2}>\sum_{k=t}^{t+j-4}
  u_{t,k}-\sum_{k=0}^{j-5} u_{t,k}+u_{1,t-1}$ for $5\leq j\leq c$.
\end{itemize} 
Then a general hyperplane section of $S$ is good determinantal 
with degree matrix ${\mathcal U}$.
\end{prop}

\begin{proof}
Let $H$ be a hyperplane that meets $S$ properly and let $C=S\cap
H$. Let $D$ be a general section of $S$. Then we have a flat family
of subschemes $D_s\subseteq\PP^n$ such that for all $s$ $D_s$ is a 
section of $S$ by a hyperplane that meets it properly, $D_0=C$ and $D_1=D$. 
Consider the Hilbert scheme $Hilb^p(\PP^n)$, where $p$ is the Hilbert polynomial of
$C$. Under our assumptions, Proposition~10.7 in~\cite{kl01} and the
results in Section~4 of~\cite{kl05}, we have that 
$dim_C Hilb^p(\PP^n)=dim W$, where $W\subseteq
Hilb^p(\PP^n)$ is the locus of good determinantal schemes whose
degree matrix is the same as the one of $C$. Moreover, $W$ is
irreducible, therefore its closure is an irreducible component of
$Hilb^p(\PP^n)$. Since $C$ is a smooth
point of $Hilb^p(\PP^n)$, we have that the irreducible component of
$Hilb^p(\PP^n)$ containing $C$ contains $D$ as well. Since $C$ belongs
to the interior of $W$, then $D_s$ belongs to $W$ for a generic value
of $s$. Therefore a general hyperplane section of $S$ is good determinantal 
with degree matrix ${\mathcal U}$.
\end{proof}

\begin{ex}
Let $S\subseteq\PP^7$ be a fourfold and let $H$ be a hyperplane that
meets $S$ properly. Let $C\subseteq H=\PP^6$ be a threefold whose 
saturated ideal is generated by the maximal minors of a generic matrix 
of linear forms of size $3\times 5$. $C$ is a smooth scroll over
$\PP^2$, and it has Hilbert polynomial
$p(t)=\frac{5}{3}t^3+4t^2+\frac{10}{3}t+1$. It follows from
Proposition~5.4 of~\cite{be05} that the Hilbert scheme $Hilb^p(\PP^6)$
has an irreducible component $\mathcal H$ of dimension 72, whose
general element is good determinantal and defined by the maximal
minors of a $3\times 5$ matrix of linear forms. $C$ is unobstructed
and it belongs to the interior of the locus of good determinantal
schemes as above (whose closure is $\mathcal H$). Then by
Proposition~\ref{sect} a general hyperplane section of $S$ is good
determinantal and defined by the maximal
minors of a $3\times 5$ matrix of linear forms.
\end{ex}

We saw that a scheme with good determinantal general hyperplane
section does not need to be good determinantal. However, it is good
determinantal up to flat deformation.

\begin{thm}\label{flatfam}
Let $S\subseteq\PP^n$, be an aCM scheme and let $C$ be a proper 
hyperplane section of $S$. Then one can find a flat family $T_s$ whose
elements all have $C$ as a proper hyperplane section, and such that
$T_1=S$ and $T_0$ is a cone over $C$.
In particular, if $C$ is standard (resp. good) determinantal, one can 
find a flat family $T_s$ whose
elements all have $C$ as a proper hyperplane section, and such that
$T_1=S$ and $T_0$ is standard (resp. good) determinantal. 
\end{thm}

\begin{proof}
By assumption $C=S\cap H$ for some hyperplane $H$ that meets $S$
properly. With no loss of generality, we can assume that $H$ is the
hyperplane of equation $x_{n+1}=0$. Let
$C\subseteq \PP^n=H\subseteq\PP^{n+1}$. Let $C'$ be the cone over $C$,
so that $H$ intersects $C'$ properly and $C'\cap H=C$. 
Then if $I_S$ has a minimal system of generators $F_1,\ldots,F_m$, 
then $I_{C'}$ has $F_1(x_0,\ldots,x_n,0),\ldots,
F_m(x_0,\ldots,x_n,0)$ as a minimal system of generators. 
Consider the flat family $T_s$ of schemes with homogeneous saturated
ideal $$I_{T_s}=(F_1(x_0,\ldots,x_n,sx_{n+1}),\ldots,
F_m(x_0,\ldots,x_n,sx_{n+1})).$$
Then $S_0=C'$ and $S_1=S$.
The graded Betti numbers are constant in the family, since the graded
Betti numbers of $C'$ and $S$ coincide by assumption, and for $s\neq
0$ $T_s$ and $S$ only differ by a change of coordinates. 
Moreover, $T_s\cap H=C$ since for all $s$
$$I_{T_s}+(x_{n+1})/(x_{n+1})=$$
$$(F_1(x_0,\ldots,x_n,0),\ldots,
F_m(x_0,\ldots,x_n,0),x_{n+1})/(x_{n+1})=I_{C|H}.$$
\end{proof}

In Theorem~\ref{flatfam} we cannot conclude that $S$ belongs to the
closure of the locus of the Hilbert scheme consisting of good
determinantal schemes.
This is connected to the fact that we cannot prove that a general element
of the flat family that we construct is good determinantal. Indeed this is not
necessarily the case, as the next example shows.

\begin{ex}
Let $V\subseteq\PP^5$ be the Veronese variety, let $C'\subseteq\PP^5$
be a cone over a rational normal quartic curve in $\PP^4$. Let 
$$M_s=\left(\begin{array}{ccc} 
x_0 & x_1 & x_2 \\
x_1 & (1-s)x_2+sx_3 & x_4 \\
x_2 & x_4 & x_5
\end{array}\right)$$
and let $T_s$ be the surface in $\PP^5$ with saturated ideal
$I_{T_s}=I_2(M_s).$
Then $S_0=C'$ while $T_s\cong V$ for $s\neq 0$. So the general element
of the flat family is not standard determinantal. Moreover, a
dimension count shows that a generic good determinantal scheme belongs
to a different component of the Hilbert scheme from the one containing
$V$. In fact, the dimension of the Hilbert scheme at $V$ is 27, while
the dimension of the component which is the closure of the locus of
good determinantal schemes is 29 (the latter can be computed using the
formulas in~\cite{kl05}). In particular $C'$ is not unobstructed
(notice that unobstructedness results such as Corollary~10.15
in~\cite{kl01} do not apply to this setting, since $C'$ is not a
Cartier divisor of the scheme defined by the matrix obtained by 
deleting a column of $M_0$). 
Notice moreover that a general hyperplane
section of $T_s$ is a rational quartic curve in $\PP^4$ for all $s$. 
\end{ex}

\section{The determinantal property via basic 
double linkage}

In this section we show how to produce a standard or good
determinantal scheme by basic double link from another determinantal
scheme. We also show how to produce a non standard determinantal
scheme by basic double link from a non standard determinantal
scheme. Putting these together, one can start from a scheme which is non
standard determinantal and whose general hyperplane section is
standard determinantal, and produce another scheme with the same
property. We refer the reader to Proposition~5.4.6 in~\cite{mi98b} for
the definition and facts about basic double links.

\begin{thm}\label{det}
Let $C\subseteq S\subseteq\PP^n$ be standard determinantal schemes, 
such that $C$ has codimension $1$ in $S$. 
Assume that for a suitable choice of defining matrices $M$ and $N$ for $C$ and
$S$, either $M$ is obtained from $N$ by deleting a row or $N$ is
obtained from $M$ by deleting a column. Then
a basic double link $D$ of $C$ on $S$ is standard determinantal. 
Moreover, if $C$ is good 
determinantal then a basic double link $D$ of $C$ on $S$ via a generic
hypersurface is good determinantal. In this sense, the property of being
standard/good determinantal is preserved under basic double linkage.
\end{thm}

\begin{proof}
Let $C\subseteq S\subseteq\PP^n$ be standard (resp. good)
determinantal schemes, where the saturated ideal of $C$ is generated
by the maximal minors of a $t\times (t+c)$ matrix $M$.  $I_C=I_t(M)$ and 
$C$ is standard determinantal, i.e. it has codimension $c+1$.
Assume that the matrix $N$
defining $S$ is obtained from the one of $C$ by adding a row, 
$I_S=I_{t+1}(N)$. $S$ has 
codimension $c$ by assumption. Notice that $S$ is good determinantal by
construction, in particular it is generically complete intersection (see \cite{kr00}, 
Remark 3.5). $I_{t+1}(N)\subseteq I_t(M)$, so
$S\supseteq C$. Let $D$ be
a basic double link of $C$ on $S$, $D=C\cup (S\cap F)$ for some
hypersurface $F$ that meets $S$ properly. If $C$ is good
determinantal, then after applying generic invertible
row operations to $M$ we have a submatrix $M'\subseteq M$
whose maximal minors define a standard determinantal scheme $U$. $M'$ is
obtained from $M$ by deleting a row. If we apply the same row
operations to $N$ and delete the corresponding row, we obtain
$N'\subseteq N$. The ideal of maximal minors of $N'$ defines a scheme
$V$ which is standard determinantal of codimension $c+1$ (as
$N$ is the defining matrix of a good determinantal scheme). We assume
that $F$ meets $V$ properly as well. Notice that this holds for a
generic choice of $F$. The saturated ideal of $D$ is then
$$I_D=I_S+F\cdot I_C$$ (see Proposition 5.4.5 in \cite{mi98b}), so it is
minimally generated by the maximal minors of the matrix obtained by adding to
$N$ a column vector, whose entries are all equal to $0$, except for an
entry equal to $F$. In other words, let $M=(m_{i,j})_{i=1,\ldots,t;\; j=1,\ldots,t+c}$ and 
$N=(n_{i,j})_{i=1,\ldots,t;\; j=1,\ldots,t+c}$, with $n_{i,j}=m_{i,j}$ for $i\leq k-1$, 
$n_{i,j}=m_{i-1,j}$ for $i\geq k+1$ (inserting a row in position $k$). If
$deg(n_{k,l-1})\leq deg(F)\leq deg(n_{k,l})$, then the defining matrix of $D$ is
$O=(o_{i,j})$ with $o_{i,j}=n_{i,j}$ for $j\leq l$, $o_{k,l}=F$, $o_{i,l}=0$ for 
$i\neq k$ and $o_{i,j}=n_{i,j-1}$ for $j\geq l+1$. Therefore $M\subseteq
N\subseteq O$, $N$ is obtained from $O$ by removing a column and $M$ is
obtained from $N$ by removing a row.
If $C$ is good determinantal, then we have $M'\subseteq M$
whose maximal minors define the standard determinantal scheme $U$. $M'$ is
obtained from $M$ by deleting a row. If we apply the same row and
column operations to $O$ and delete the corresponding row, we obtain
$O'\subseteq O$. The ideal of maximal minors of $O'$ defines a scheme
which is a basic double link of $U$ on $V$. In fact
$$I_t(O')=I_t(N')+F\cdot I_{t-1}(M')=I_V+F\cdot I_U.$$
Recall that by assumption $F$ meets $V$ properly. In particular,
$I_t(O')$ defines a standard determinantal scheme of codimension
$c+2$. This proves that $D$ is good determinantal.

Assume now the matrix $N$ that defines $S$ is
obtained from $M$ by deleting the $k$-th column.
$I_S=I_t(N)$ and $S$ has codimension $c$ by assumption. 
Notice that all the minimal generators of $I_S$ are
also minimal generators of $I_C$. Moreover, $S$ is good determinantal
(as shown in \cite{kl01}, Theorem~3.6).
$I_t(N)\subseteq I_t(M)$, so $S\supseteq C$. Let $D$ be a basic double link of $C$ on
$S$, $D=C\cup (S\cap F)$ for some hypersurface $F$ that meets $S$
properly. The saturated ideal of $D$ is $I_D=I_S+F\cdot I_C$ (see
Proposition~5.4.5 in~\cite{mi98b}), so it is minimally generated by
the maximal minors of the matrix $O$ obtained by adding to $N$ the $k$-th
column of $M$, after multiplying all of the entries by $F$. 
If $C$ is good determinantal, then after applying generic invertible
row operations to $M$ we have a submatrix $M'\subseteq M$
whose maximal minors define a standard determinantal scheme $U$. $M'$ is
obtained from $M$ by deleting a row. If we apply the same row 
operations to the matrix $O$ and delete the corresponding row, we obtain
$O'\subseteq O$. The ideal of maximal minors of $O'$ defines a scheme
which is a basic double link of $U$, in particular it has codimension
$c$ hence it is standard determinantal. Therefore $D$ is good determinantal.  
\end{proof}

We now give an example of how
one can systematically produce families of schemes 
which are not standard determinantal. This can be achieved 
by taking a basic double link of a scheme $C$ which 
is not standard determinantal on a standard determinantal scheme
$S$. Of course one needs to check that the result is not standard
determinantal, since clearly Theorem~\ref{det} does not guarantee
it. Let $H$ be a hyperplane that meets $C,D$ and $S$ properly.
In order to guarantee that the basic double link $D\cap H$ of $C\cap
H$ on $S\cap H$ is standard determinantal, we can lift a basic double
link of the type described in Theorem~\ref{det}
from $C\cap H$ to $C$.

\begin{ex}
Let $V\subseteq\PP^5$ be the Veronese surface
$$I_V=I_2\left(\begin{array}{ccc} 
x_0 & x_1 & x_2 \\
x_1 & x_5 & x_3 \\
x_2 & x_3 & x_4 
\end{array}\right)\subseteq k[x_0,\ldots,x_5].$$
Let $S\subseteq\PP^5$ be the threefold defined by 
$$I_S=I_3\left(\begin{array}{cccc} 
x_0 & x_1 & x_2 & x_3 \\
x_1 & x_5 & x_3 & x_4 \\
x_2 & x_3 & x_4 & x_0
\end{array}\right).$$
Then $S$ is good determinantal and contains $V$. 
Let $F$ be a general linear form. Then a basic double
link $W=V\cup(S\cap F)$ of $V$ on $S$ is not standard
determinantal. This can be checked by computing the cardinality of a
minimal system of generators of $I_W^2$ and counting Pl\"ucker
relations (as done in Example~\ref{verodeform}).

Let $C\subseteq\PP^4$ be a smooth rational normal curve
$$I_C=I_2\left(\begin{array}{cccc} 
x_0 & x_1 & x_2 & x_3 \\
x_1 & x_2 & x_3 & x_4
\end{array}\right)\subseteq k[x_0,\ldots,x_4].$$
Let $H\subseteq\PP^5$ by the hyperplane of equation $x_2-x_5=0$. Then
$H$ meets $V$ properly and $C=V\cap H\cong\PP^4$. Moreover, if
$T=S\cap H$, then $D=W\cap H=C\cup(T\cap F)$ is a basic double link
of $C$ on $T$. $D$ is good determinantal by Theorem~\ref{det}. The
saturated ideal of $D$ is 
$$I_D=I_3\left(\begin{array}{ccccc} 
x_0 & x_1 & x_2 & x_3 & 0 \\
x_1 & x_2 & x_3 & x_4 & 0 \\
x_2 & x_3 & x_4 & x_0 & \overline{F}
\end{array}\right)$$
where we denote by $\overline{F}$ the equation of $F$ restricted to
$H$.
\end{ex}

Next we show in an example how one can use a similar construction to
produce a scheme that is not standard determinantal and 
whose general hyperplane section is good determinantal.

\begin{ex}\label{gensectbdl}
Consider the curve $C\subseteq\PP^{n+1}$ of
Proposition~\ref{n+1curve}. We use the same notation as in the proof
of the proposition. We saw that 
$$I_C=x_0(x_2,\ldots,x_{n+1})+\sum_{1\leq i<j\leq
  n}(x_ix_j).$$ Let $S\supseteq C$ be the surface cut out by all the
squarefree monomials of degree 2 in $x_0,x_2,\ldots,x_n$. $S$ is
good determinantal by Lemma~\ref{sqfr}. Let $L$ be a
hyperplane that meets $S$ properly, 
let $D=C\cup(S\cap L)$. To simplify the computation, we let $L=x_1$. $D$ is a basic double
link of $C$ on $S$, and has saturated ideal $$I_D=x_1I_C+I_S=(x_0x_1x_{n+1})
+x_1^2(x_2,\ldots,x_n)+\sum_{i,j\in\{0,2,\ldots,n\},\; i<j}(x_ix_j).$$
We now sketch the proof that $D$ is not standard determinantal.
In order to show it, we proceed as in the proof of
Proposition~\ref{n+1curve} and examine the last matrix in a minimal free
resolution of the ideal of $D$. 
We can follows the same steps as in the proof of
Proposition~\ref{n+1curve}, taking into account the fact that the minimal generators
$x_0x_{n+1}$ and $x_1(x_2,\ldots,x_n)$ are replaced by their multiples by $x_1$.
Therefore we can write the last matrix in a minimal free resolution of $I_D$ as 
\begin{equation}\label{lastD}\left(\begin{array}{cccccc}
x_2 & 0 & \ldots & \ldots & \ldots & 0 \\
-x_3 & -x_1^2 & 0 & \ldots & \ldots & 0 \\
\vdots & 0 & \ddots & & & \vdots \\
\vdots & \vdots & \ddots &\ddots & & \vdots \\
(-1)^n x_n & 0 & \ldots & 0 & -x_1^2 & 0 \\
(-1)^{n+1} x_1x_{n+1} & 0 & \ldots & 0 & 0 & -x_1^2 \\
0 & x_0 & 0 & \ldots  &\ldots & 0 \\
\vdots & 0 & \ddots & & & \vdots \\
\vdots & \vdots & & &\ddots & 0 \\
0 & 0 & \ldots & 0 & 0 & x_0 \\
0 & & & & & \\
0 & & & M' & & \\
0 & & & & & \\
\end{array}\right).\end{equation}
The matrix $M'$ is obtained from the matrix $M$ in (\ref{matx}) by replacing 
each occurrence of $x_1$ by $x_1^2$.
Then one checks that the ideal of $2\times
2$ minors of the matrix (\ref{lastD}) is monomial, and it does not contain 
any pure power of $x_{n+1}$.
However, it contains all the monomials of degree $2$ in
$x_0,x_2,\ldots,x_n$, as well as $x_1^4$, $x_1^3x_{n+1}$ and $x_1^2x_i$ for 
$i=0,2,\ldots,n$. As in Proposition~\ref{n+1curve}, one can write
down the last matrix in a minimal free resolution of the ideal of
maximal minors of a $2\times(n+1)$ matrix of indeterminates $z_0,\ldots,z_{2n+1}$. 
The matrix has been explicitly described in (\ref{matx2}). One can check
that the ideal of $2\times 2$ minors of (\ref{matx2}) is $(z_0,\ldots,z_{2n+1})^2$. 
Therefore, we conclude that $D$ is not standard determinantal by a specialization 
argument as in Proposition~\ref{n+1curve}. If $I_D$ is the ideal of
$2\times 2$ 
minors of a $2\times (n+1)$ matrix of linear forms, then 
the entries of the matrix do not involve $x_{n+1}$, which is a contradiction.

We show that a general hyperplane section of $D$ is good determinantal.
Let $H\subseteq\PP^{n+1}$ be a general hyperplane of equation
$x_{n+1}-h$. Let $X=C\cap H$,
$Y=D\cap H$, and $E=S\cap H$. Then $X,Y$ are zero-dimensional
subschemes of $H\cong\PP^n$, $X,Y\subseteq E$. $x_0,\ldots,x_n$ are
coordinates on $H$ and 
$$I_{X|H}=x_0(x_2,\ldots,x_n,h)+\sum_{1\leq i<j\leq n}(x_ix_j)=$$
$$x_0(y,x_2,\ldots,x_n)+y(x_2,\ldots,x_n)+\sum_{2\leq i<j\leq n}(x_ix_j).$$
Here $y=\alpha x_0+\beta x_1$ for generic $\alpha,\beta\in k^*$. 
Then $I_{X|H}$ is generated by the squarefree monomials of degree 2
in $x_0,y,x_2,\ldots,x_n$, hence it corresponds to $n+1$ generic
points in $\PP^n$. So $I_{X|H}$ is the ideal of maximal minors of the matrix 
$$\left(\begin{array}{ccccc}
x_0 & y & x_2 & \ldots & x_n \\
x_0 & \gamma y & \gamma^2 x_2 & \ldots & \gamma^n x_n 
\end{array}\right)$$
for $\gamma\in k^*$ generic. 
The ideal $I_{S|H}$ is generated by the maximal minors of 
$$\left(\begin{array}{cccc}
x_0 & x_2 & \ldots & x_n \\
x_0 & \gamma^2 x_2 & \ldots & \gamma^n x_n
\end{array}\right).$$
Therefore Theorem~\ref{det} applies, and
$Y\subseteq H\cong\PP^n$ is good determinantal with defining matrix
$$\left(\begin{array}{ccccc}
x_0 & yx_1 & x_2 & \ldots & x_n \\
x_0 & \gamma yx_1 & \gamma^2 x_2 & \ldots & \gamma^n x_n 
\end{array}\right).$$
\end{ex}

\end{document}